\documentclass[12pt]{article}

\DeclareMathAlphabet{\mathpzc}{OT1}{pzc}{m}{it}

\usepackage{amsmath}
\usepackage{amssymb}
\usepackage{amsfonts}
\usepackage{amsthm}
\usepackage{mathrsfs}
\usepackage{ifsym}
\usepackage{fontenc}
\usepackage{textcomp}
\usepackage{tipa}
\usepackage{enumerate}
\usepackage[pdftex]{color, graphicx}

\begin{document}

\title{{\bf  Product formulas for the $5$-division points on the Tate normal form and the Rogers-Ramanujan continued fraction}}         
\author{Patrick Morton}        
\date{}          
\maketitle

\begin{abstract}
Explicit formulas are proved for the $5$-torsion points on the Tate normal form $E_5$ of an elliptic curve having $(X,Y)=(0,0)$ as a point of order $5$.  These formulas express the coordinates of points in $E_5[5] - \langle(0,0)\rangle$ as products of linear fractional quantities in terms of fifth roots of unity and a parameter $u$, where the parameter $b$ which defines the curve $E_5$ is given as $b=(\varepsilon^5 u^5- \varepsilon^{-5})/(u^5+1)$ and $\varepsilon = (-1+\sqrt{5})/2$.  If $r(\tau)$ is the Rogers-Ramanujan continued fraction and $b=r^5(\tau)$, then the coordinates of points of order $5$ in $E_5[5] - \langle(0,0)\rangle$ are shown to be products of linear fractional expressions in $r(5\tau)$ with coefficients in $\mathbb{Q}(\zeta_5)$.
\end{abstract}

\section{Introduction.}

In previous papers, several new formulas for the $3$-division points on the Deuring normal form
$$E_3: \ Y^2+\alpha XY +Y = X^3,$$
and the $4$-division points on the Tate normal form
$$E_4: \ Y^2+XY + bY = X^3+bX^2,$$
have recently been given.  For the curve $E_3$, the point
$$\left( X, Y \right) = \left(\frac{-3\beta}{\alpha(\beta-3)},\frac{\beta-3\omega}{\beta-3}\right)$$
represents the six points of order $3$ in $E_3[3]-\langle (0,0) \rangle$, where $\omega$ is one of the two primitive cube roots of unity and $(\alpha, \beta)$ lies on the Fermat cubic
$$Fer_3: \ 27X^3+27Y^3=X^3Y^3.$$
(See \cite{mor3}.)  Setting $b=1/\alpha^4$ in the equation for $E_4$, a point of order $4$ in $E_4[4]-\langle (0,0) \rangle$ is the point
$$
(X,Y) = \left(-\beta_1 \beta_2 \beta_3, \beta_1^2 \beta_2^2 \beta_3\right), \eqno{(1.1)}
$$
where
$$\beta_n=\frac{\beta+2i^n}{2\beta}, \ \ i = \sqrt{-1},$$
and the point $(\alpha, \beta)$ lies on the Fermat quartic
$$Fer_4: \ 16X^4+16Y^4=X^4 Y^4.$$
Replacing $\beta$ by $i \beta$ (so $\beta_n$ becomes $\beta_{n-1}$) and $i$ by $-i$ in the above formula yields $8$ of the $12$ points of order $4$ in $E_4[4]$.  The other points of order $4$ are $(0,0), (0,-b) \in \langle (0,0) \rangle$ and the two points $(-2b,2\beta_1\beta_3b)$ and $(-2b,2\beta_2\beta_4 b)$.
(See \cite{ly}, \cite{lym}.) \medskip

Similar formulas have been given for the $6$-torsion points on the Tate normal form $E_6$, by Lynch \cite{ly}. This normal form is
$$E_6: \ Y^2+aXY+bY=X^3+bX^2, \ \ b=-(a-1)(a-2).$$
Lynch's formulas express the coordinates of $6$-torsion points on $E_6$ as products of linear fractional quantities in $\alpha$ and $\beta$ and a cube root of unity $\omega$, where $(\alpha,\beta)$ is a point on the elliptic curve
$$Y^2=X^3+1$$
and the parameter $a$ is given by
$$a=\frac{10\beta^2-18}{9(\beta^2-1)}=\frac{10\alpha^3-8}{9\alpha^3}.$$
The exact formulas are somewhat complicated; these can be found in \cite{ly}. \medskip

In this note I will prove analogous formulas for the non-trivial points of order $5$ on the Tate normal form
$$E_5(b) : \ Y^2+(1+b)XY+bY=X^3+bX^2,$$
on which the point $(0,0)$ is a point of order $5$.  (See the discussion in \cite{mor} for more on the Tate normal form.)  These formulas are similar to the expressions (1.1) for the points of order $4$ on the curve $E_4$, in that they express the $X$ and $Y$ coordinates of points in $E_5(b)[5]-\langle(0,0)\rangle$  as products of linear fractional quantities in a parameter $u$, where
$$b=\frac{\varepsilon^5u^5+\bar \varepsilon^5}{u^5+1}, \ \ \varepsilon = \frac{-1+\sqrt{5}}{2}, \ \ \bar \varepsilon = \frac{-1-\sqrt{5}}{2};$$
and the coefficients in these linear fractional expressions lie in the field $\mathbb{Q}(\zeta_5)$ of fifth roots of unity.  These expressions are quite a bit simpler than the formulas given by Verdure \cite{v}, in which the $Y$-coordinates of the $5$-torsion points are expressed in terms of a formal root $x_0$ of the $5$-division polynomial.  In this paper, the quantity
$$u^5=-\frac{b-\bar \varepsilon^5}{b - \varepsilon^5},$$
which is up to sign the same as Verdure's Kummer element (see Theorem 5 in \cite{v}), arises naturally in the process of solving the quintic equation $g(X)=0$ below using Watson's method (see \cite{lsw} and Section 5 of this paper).  \medskip

The expressions given in Theorem 2.1 below also allow one to check ``by hand" that these points do indeed have order 5, assuming that they represent points on $E_5$.  (See Theorem 3.1 and the discussion in Section 3.)  These formulas will be used in a forthcoming paper to prove the case $p=5$ of the conjectures stated in \cite{mor4} and \cite{mor5}.
\medskip

These formulas also have a strong connection to the Rogers-Ramanujan continued fraction, which is
$$r(\tau) = \frac{q^{1/5}}{1+\frac{q}{1+\frac{q^2}{1+ \frac{q^3}{1+\cdots}}}}=\frac{q^{1/5}}{1+} \ \frac{q}{1+} \ \frac{q^2}{1+} \ \frac{q^3}{1+} \dots,$$
and whose value is the modular function for $\Gamma(5)$ given by
$$
r(\tau) = q^{1/5} \prod_{n \ge 1}{(1-q^n)^{(n/5)}},  \ \ q = e^{2 \pi i \tau}, \ \ \tau \in \mathbb{H}.\eqno{(1.2)}
$$
(The symbol $\left(\frac{n}{5}\right)$ in the exponent is the Legendre symbol and $\mathbb{H}$ is the upper half-plane.  See \cite{anb}, \cite{ber}, and \cite{duke} and the references in the latter paper.). From the formulas of \cite{duke} it follows easily that if $b=r^5(\tau)$, then the parameter $u$ described above may be taken to be
$$u =\frac{1}{\varepsilon r\left(\frac{-1}{5\tau}\right)}=-\frac{r(5\tau)-\bar \varepsilon}{r(5\tau)-\varepsilon}.\eqno{(1.3)}$$
This yields the following. \bigskip

\noindent {\bf Theorem 1.1.} {\it If $b=r^5(\tau)$, for $\tau \in \mathbb{H}$, then 
$$X=\frac{-\varepsilon}{\sqrt{5}} \ \frac{r^4(5\tau)-3r^3(5\tau)+4r^2(5\tau)-2r(5\tau)+1}{r^2(5\tau)+r(5\tau)+\varepsilon^2}\eqno{(1.4)}$$
is the $X$-coordinate of a point $P=(X,Y)$ of order $5$ on the elliptic curve $E_5(b)$, which is not in the group $\langle (0,0) \rangle$.  The $Y$-coordinates of $P$ and $-P$ are products of linear fractional expressions in $r(5\tau)$ with coefficients in $\mathbb{Q}(\zeta_5)$; and the same holds for the coordinates of all points in $E_5(b)[5]-\langle(0,0)\rangle$}.\bigskip

This formula (1.4) is closely related to a well-known identity of Ramanujan:
$$\frac{r^5(\tau)}{r(5\tau)}=\frac{r^4(5\tau)-3r^3(5\tau)+4r^2(5\tau)-2r(5\tau)+1}{r^4(5\tau)+2r^3(5\tau)+4r^2(5\tau)+3r(5\tau)+1}.$$
(See \cite{ber}, p. 167; also \cite{duke}, equation (7.4), except that the term $r(5\tau)$ on the left side of (7.4) should be $r(\tau)$.)  This identity allows us to express the formula for $X$ in the following form:
$$X=\frac{-\varepsilon}{\sqrt{5}} \ \frac{r^5(\tau)}{r(5\tau)} (r^2(5\tau)+r(5\tau)+\bar \varepsilon^2).$$
One of the corresponding $Y$-coordinates is given by the formula
$$Y_1=\left(\frac{\zeta-1}{\sqrt{5}}\right)^3 \frac{r^5(\tau)}{r(5\tau)} (r^4(5\tau)-3r^3(5\tau)+4r^2(5\tau)-2r(5\tau)+1) \cdot Z,$$
where
$$Z=\frac{r(5\tau)-(\zeta^3+\zeta^4)}{(r(5\tau)-(\zeta^2+\zeta^4))(r(5\tau)-(1+\zeta^3))},\ \ \ \zeta=\zeta_5;$$
and the other is obtained by replacing $\zeta$ in this formula by $\zeta^4$.  (See equations (4.3) and (4.4) below.)  Replacing $\zeta$ by $\zeta^2$ and interchanging $\varepsilon$ and $\bar \varepsilon$ in these formulas yields the coordinates of the points $\pm 2P$.  Thus the subgroup $\langle P \rangle$ is determined by the action of $\textrm{Gal}(\mathbb{Q}(\zeta)/\mathbb{Q})$ on $P$.  (See Theorem 3.1.)
\medskip

In the related paper \cite{mor6} this connection with $r(\tau)$ will be applied to show the following result in the theory of complex multiplication.  If $-d=d_K f^2$ is the discriminant of the order $\textsf{R}_{-d}$ of conductor $f$ in the quadratic field $K=\mathbb{Q}(\sqrt{-d})$, where $\left(\frac{-d}{5}\right)=1$ ($d \neq 4f^2$); and if $\tau$ has the value
$$\tau=\frac{v+\sqrt{-d}}{10}, \ \ v^2+d \equiv 0 \ (\textrm{mod} \ 4 \cdot 5^2), \ \ (v,f) = 1;$$
then the unit $r(5\tau) = r\left(\frac{v+\sqrt{-d}}{2}\right)$ generates the field $F=\Sigma_5 \Omega_f$ over $\mathbb{Q}$, where $\Sigma_5$ is the ray class field of conductor $5$ over $K$ and $\Omega_f=K(j(\tau))$ is the ring class field of conductor $f$ over $K$.  Furthermore, for some primitive $5$-th root of unity $\zeta$,
$$\mathbb{Q}(r(\tau))=\Sigma_{\wp_5'}\Omega_f, \ \ \mathbb{Q}(\zeta r\left(\frac{-1}{5\tau}\right))=\Sigma_{\wp_5}\Omega_f, \ \ (\zeta \neq 1),$$
where $\wp_5$ is the prime ideal divisor of $5$ in $K$ for which $\wp_5 \mid 5\tau$ and $\wp_5'$ is the conjugate prime ideal.  In particular, $\mathbb{Q}(r(5\tau))$ is a normal extension of $\mathbb{Q}$, while $\mathbb{Q}(r(\tau))$ is not normal over $\mathbb{Q}$, though both are abelian over $K$.  At any rate, values of the Rogers-Ramanujan function $r(\tau)$ turn out to yield generators of small height for class fields of quadratic fields $K$ in which the prime $5$ splits.  The reader is referred to \cite{mor6} for a list of the minimal polynomials of these values for small values of $d$.

\section{Points of order 5 on $E_5(b)$.}

The $X$-coordinates of points of order 5 on $E_5(b)$ which are not in the group 
$$\langle (0,0) \rangle = \{O, (0,0), (0,-b), (-b,0), (-b,b^2)\}$$
are roots of the polynomial
\begin{align*}
D_5(x)&=5x^{10}+(5+25b+5b^2)x^9+(1+38b+44b^2+7b^3+b^4)x^8\\
&+(9b+127b^2+26b^3+3b^4-b^5)x^7+(36b^2+248b^3+19b^4-3b^5+b^6)x^6\\
&+(84b^3+322b^4+71b^5+3b^6-b^7)x^5+(126b^4+293b^5+94b^6+12b^7+b^8)x^4\\
&+(125b^5+180b^6+50b^7+5b^8)x^3+(80b^6+65b^7+10b^8)x^2\\
&+(30b^7+10b^8)x+5b^8.
\end{align*}
This follows easily from \cite{si1}, Exercise 3.7 (p. 105), applied to the curve $E_5$, after factoring out $x(x+b)$ from the polynomial $\psi_5(x)$.  (But note that the formula for $b_2$ on p. 42 should be $b_2=a_1^2+4a_2$.)  This polynomial factors into $5$ times the product of two polynomials
\begin{align*}
g(X)&=X^5+\frac{1}{20}(\alpha-5)(-3-\alpha-7b+3b\alpha -2b^2)X^4\\
&+\frac{\alpha}{5}b(1+2\alpha-11b+4b\alpha-b^2)X^3\\
&+\frac{1}{10}(\alpha-5)b^2(-9-2\alpha-6b+b\alpha-b^2)X^2+(3b^3+b^4)X+b^4,
\end{align*}
where $\alpha^2=5$.  Using Watson's method of solving the for the roots of a quintic equation from \cite{lsw}, we find that the roots of $g(X)$ are given by
\begin{align*}
X&=\frac{(5-\alpha)}{100} \{(-18+8\alpha-12b+6b\alpha-2b^2)u^4+(-7+3\alpha+12b-4b\alpha+2b^2)u^3\\
&+(-3+\alpha-7b+7b\alpha-2b^2)u^2+(-2+22b+2b^2)u-3-\alpha-7b+3b\alpha-2b^2 \}\\
&=\frac{(5-\alpha)}{100}(A_4u^4+A_3u^3+A_2u^2+A_1u+A_0),
\end{align*}
where
\begin{equation*}
u^5=\phi(b) = \frac{2b+11+5\alpha}{-2b-11+5\alpha}=\frac{b-\bar \varepsilon^5}{-b+\varepsilon^5},
\end{equation*}
$$\varepsilon = \frac{-1+\alpha}{2}=\zeta+\zeta^4, \ \ \bar \varepsilon = \frac{-1-\alpha}{2}=\zeta^2+\zeta^3,$$
and $\zeta$ is a primitive $5$-th root of unity.  (The details of Watson's method applied to the polynomial $g(X)$ are given in the appendix.  Note that $\varepsilon$ and $\bar \varepsilon$ are the quadratic Gaussian periods for $\mathbb{Q}(\zeta)$.)  This may be verified on Maple by plugging the expression for $X$ into $g(X)$, and using the formula
$$b = \frac{\varepsilon^5u^5+\bar \varepsilon^5}{u^5+1}\eqno{(2.1)}$$
for $b$ in terms of $u$.  \medskip

Using this formula for $b$, the above value of $X$ can also be written as
$$X=\frac{(-7+3\alpha)}{4}\frac{(-2u^2+(1+\alpha)u-3\alpha-7)(2u^2+(4+2\alpha)u+3\alpha+7)}{(-2u^2+(1+\alpha)u-2)(u+1)^2}.\eqno{(2.2)}$$
The formulas (2.1) and (2.2) show that there are $10$ such values, since replacing $u$ by $\zeta^i u$ (and leaving $\alpha$ unchanged), or replacing $\alpha$ by $-\alpha$ and $u$ by $1/(\zeta^i u)$ gives the other $X$-coordinates.  It is easy to see that these transformations yield distinct points in $E_5(b)[5]$, since the $X$-coordinates have distinct sets of poles.  Setting $\alpha=\zeta-\zeta^2-\zeta^3+\zeta^4$, this expression factors:
$$X=-\varepsilon^4 \frac{[u-(1+\zeta)^2][u-\zeta(1+\zeta)^2][u-\zeta^2(1+\zeta)^2][u-\zeta^3(1+\zeta)^2]}{(u+\zeta^2)(u+\zeta^3)(u+1)^2}.$$
The zeros and poles of this function of $u$ are all units in $\mathbb{Q}(\zeta_5)$.  We will now show that the corresponding $Y$-coordinates factor in a similar way.  We derive the following theorem using calculations in an extension of the field $\mathbb{Q}(\zeta, b)$, but the formulas themselves are valid over any field whose characteristic is different from $5$.
\bigskip

\noindent {\bf Theorem 2.1.}
{\it If $b = \frac{\varepsilon^5u^5+\bar \varepsilon^5}{u^5+1}$, the $X$-coordinates of the points of order $5$ in $E_5[5]-\langle (0,0) \rangle$ are given by the formula}
$$X=-\varepsilon^4 \frac{[u-(1+\zeta)^2][u-\zeta(1+\zeta)^2][u-\zeta^2(1+\zeta)^2][u-\zeta^3(1+\zeta)^2]}{(u+\zeta^2)(u+\zeta^3)(u+1)^2},\eqno{(2.3)}$$
{\it where $\varepsilon = \frac{-1+\alpha}{2}$, $\alpha=\pm \sqrt{5}=\zeta-\zeta^2-\zeta^3+\zeta^4$, and $\zeta=\zeta_5$ is a primitive 5-th root of unity.  The corresponding $Y$-coordinates are given by}
$$Y_1=\varepsilon^7 \frac{[u-(1+\zeta)^2]^2[u-\zeta(1+\zeta)^2]^2[u-\zeta^2(1+\zeta)^2]^2[u-\zeta^3(1+\zeta)^2]}{(u+\zeta^2)^2(u+\zeta^3)(u+\zeta^4)(u+1)^3},$$
{\it and}
$$Y_2=\varepsilon^7 \frac{[u-(1+\zeta)^2][u-\zeta(1+\zeta)^2]^2[u-\zeta^2(1+\zeta)^2]^2[u-\zeta^3(1+\zeta)^2]^2}{(u+\zeta)(u+\zeta^2)(u+\zeta^3)^2(u+1)^3}.$$ \smallskip

\noindent {\it Proof.}
Putting (2.1) and (2.2) into the equation for $E_5$ yields the following equation for $Y$:
$$AY^2+BY+C=0,$$
where
\begin{align*}
A=&\frac{1}{8} (2u^2+(-1+\alpha)u+2)(-2u^2+(1+\alpha)u-2)^3(u+1)^6;\\
B=&\frac{-\varepsilon^7}{16}(-4u^3+(3+\alpha)u^2-2(1+\alpha)u+6+2\alpha) (-2u^2+(1+\alpha)u-7-3\alpha)\\
\times & (-2u^2+(1+\alpha)u-2)(2u^2+(4+2\alpha)u+7+3\alpha)^2(u+1)^3;\\
C=& \frac{\varepsilon^{14}}{64} (-2u^2+(1+\alpha)u-7-3\alpha)^3(2u^2+(4+2\alpha)u+7+3\alpha)^4.
\end{align*}
The discriminant of the quadratic is
\begin{align*}
D=&\frac{-5 \alpha \varepsilon^{13}}{64}u^4(-2u^2+(1+\alpha)u-7-3\alpha)^2(-2u^2+(1+\alpha)u-2)^2\\
\times & (2u^2+(4+2\alpha)u+7+3\alpha)^4(u+1)^6,
\end{align*}
which is $-\alpha \varepsilon=(\zeta^2-\zeta^3)^2$ times a square.  Thus, the roots of the quadratic are
$$Y=\frac{-B \pm (\zeta^2-\zeta^3) \alpha \frac{\varepsilon^6}{8}S}{2A}=\frac{-16B \pm 2 (\zeta^2-\zeta^3) (\zeta-\zeta^2-\zeta^3+\zeta^4) \varepsilon^6 S}{32A},$$
where
$$S=u^2(-2u^2+(1+\alpha)u-7-3\alpha)(-2u^2+(1+\alpha)u-2)(2u^2+(4+2\alpha)u+7+3\alpha)^2(u+1)^3.$$
Now, $1/\varepsilon=-\bar \varepsilon = -(\zeta^2+\zeta^3)$, which gives that
$$(\zeta^2-\zeta^3) (\zeta-\zeta^2-\zeta^3+\zeta^4) \varepsilon^6=\varepsilon^7 (-\zeta^3+3\zeta^2+2\zeta+1).$$
The numerator in the expression for $Y$ then becomes
\begin{align*}
-16B \ \pm \ & 2(-\zeta^3+3\zeta^2+2\zeta+1) \varepsilon^7S=\varepsilon^7 (-2u^2+(1+\alpha)u-7-3\alpha) \\
\times & (-2u^2+(1+\alpha)u-2)(2u^2+(4+2\alpha)u+7+3\alpha)^2(u+1)^3\\
\times & \{(-4u^3+(3+\alpha)u^2-2(1+\alpha)u+6+2\alpha)\pm 2(-\zeta^3+3\zeta^2+2\zeta+1)u^2\}.
\end{align*}
The quantities inside the brackets are, respectively,
\begin{align*}
(-4u^3+(3+\alpha)u^2 & -2(1+\alpha)u+6+2\alpha) + 2(-\zeta^3+3\zeta^2+2\zeta+1)u^2\\
&=-4(u+\zeta)(u+\zeta^3)(u-(1+\zeta)^2),
\end{align*}
and
\begin{align*}
(-4u^3+(3+\alpha)u^2 & -2(1+\alpha)u+6+2\alpha) - 2(-\zeta^3+3\zeta^2+2\zeta+1)u^2\\
&=-4(u+\zeta^2)(u+\zeta^4)(u-\zeta^3(1+\zeta)^2).
\end{align*}
On the other hand, the factors of the quantity $A$ are
$$2u^2+(-1+\alpha)u+2=2(u+\zeta)(u+\zeta^4),$$
while
$$-2u^2+(1+\alpha)u-2=-2(u+\zeta^2)(u+\zeta^3).$$
Now, using the factorizations
\begin{align*}
-2u^2+(1+\alpha)u-7-3\alpha & =-2(u-(1+\zeta)^2)(u-\zeta^3(1+\zeta)^2),\\
2u^2+(4+2\alpha)u+7+3\alpha & = 2(u-\zeta(1+\zeta)^2)(u-\zeta^2(1+\zeta)^2),
\end{align*}
we find the two expressions $Y_1$ and $Y_2$ stated in the theorem.  These factorizations also yield the factorization of the numerator and denominator of $X$ in (2.2).
$\square$  \bigskip

\noindent {\bf Remarks.}
The theorem shows that the quantities $X$ and $Y_i$ factor in a similar way over $\mathbb{Q}(\zeta)$ to the way that the quantity $b$ factors:
$$b=\varepsilon^5 \frac{[u-(1+\zeta)^2][u-\zeta(1+\zeta)^2][u-\zeta^2(1+\zeta)^2][u-\zeta^3(1+\zeta)^2][u-\zeta^4(1+\zeta)^2]}{(u+\zeta)(u+\zeta^2)(u+\zeta^3)(u+\zeta^4)(u+1)}.$$
The expression for $X$ may be written as
$$X=-\varepsilon^4 \frac{(u-(1+\zeta)^2)}{u+1} \frac{(u-\zeta(1+\zeta)^2)}{u+\zeta}\frac{(u-\zeta^2(1+\zeta)^2)}{u+\zeta^2}\frac{(u-\zeta^3(1+\zeta)^2)}{u+\zeta^3}\frac{u+\zeta}{u+1},$$
and the $Y_i$ may be written in a similar form.  Thus, the coordinates of $P=(X,Y_i)$ are products of linear fractional expressions in $u$. \medskip

\section{Checking the formulas.}

The curve $E_5$ is isomorphic to the curve
$$E': \ Y'^2=X^3+\frac{b^2+6b+1}{4}X^2+\frac{b(b+1)}{2}X+\frac{b^2}{4}$$
by the substitution $Y'=Y+\frac{1}{2}(1+b)X+\frac{b}{2}$.  From \cite{si1}, Ex. 3.7 the doubling formula on $E'$ is given by
$$X(2P) = \frac{X^4-(b^2+b)X^2-2b^2X-b^3}{4p(X)}, \ \ X=X(P),$$
with $p(X) =X^3+\frac{b^2+6b+1}{4}X^2+\frac{b(b+1)}{2}X+\frac{b^2}{4}$; and
$$Y'(2P) = \frac{N(X)}{16p(X)^2}Y'(P),$$
with
$$N(X)=2X^6+(b^2+6b+1)X^5+(5b^2+5b)X^4+10b^2X^3+10b^3X^2+(b^5+5b^4)X+b^5.$$

Taking the expression for $X=X(P)$ from (2.2), we have
$$X(2P)=\frac{(-7+3\alpha)}{4}\frac{(-2u^2+(1+\alpha)u-3\alpha-7)(2u^2+(4+2\alpha)u+3\alpha+7)}{(-2u^2+(1-\alpha)u-2)(u+1)^2},$$
which only differs from (2.2) in the denominator, where $\alpha$ has been replaced by its conjugate $-\alpha$.
Notice that the numerator in this formula for $X(2P)$ is
\begin{align*}
& (-7+3\alpha)(-2u^2+(1+\alpha)u-3\alpha-7)(2u^2+(4+2\alpha)u+3\alpha+7)\\
& = (28-12\alpha)u^4+(12-4\alpha)u^3+8u^2+(12+4\alpha)u+28+12\alpha.
\end{align*}
This expression is invariant (except for a factor of $u^4$) under the mapping $(\alpha \rightarrow -\alpha, u \rightarrow 1/u)$.  From (2.1) we see that this mapping also leaves the quantity $b$ invariant, and takes the denominator of $X$ (a symmetric polynomial in $u$) to the denominator of $X(2P)$ divided by $u^4$.  Hence, $X(2P)$ is the $X$-coordinate in Theorem 2.1 corresponding to the pair $(-\alpha,1/u)$, and we may state the following. \bigskip

\noindent {\bf Theorem 3.1.} {\it If $X$ is given by (2.2), the $X$-coordinate of the double of the point $P=(X,Y_i)$ on $E_5$ is obtained by applying the mapping $(\alpha \rightarrow -\alpha, u \rightarrow 1/u)$ to the expression (2.2) or $(\zeta \rightarrow \zeta^2, u \rightarrow 1/u)$ to (2.3).}
\bigskip

Since the mapping $(\alpha \rightarrow -\alpha, u \rightarrow 1/u)$ has order $2$, it is clear that $X(4P)=X(P)$ for either of the points $P=(X,Y_i)$ in Theorem 2.1.  Applying the map $\sigma = (\zeta \rightarrow \zeta^2, u \rightarrow 1/u)$ to the quantity $Y_1$ in Theorem 2.1 yields
$$Y_1^\sigma=\bar \varepsilon^7 \frac{[1-(1+\zeta^2)^2u]^2[1-\zeta^2(1+\zeta^2)^2u]^2[1-\zeta^4(1+\zeta^2)^2u]^2[1-\zeta(1+\zeta^2)^2u]}{(1+\zeta^4u)^2(1+\zeta u)(1+\zeta^3u)(u+1)^3},$$
and therefore, since $\frac{1}{(1+\zeta^2)^2}=\zeta^2(1+\zeta)^2$ and
$$(\zeta^2+\zeta^3)^7 \cdot (1+\zeta^2)^{14} \cdot \zeta = 21+13(\zeta^2+\zeta^3)=-\varepsilon^7,$$
we find that
$$Y_1^\sigma=\varepsilon^7 \frac{[u-(1+\zeta)^2]^2[u-\zeta(1+\zeta)^2][u-\zeta^2(1+\zeta)^2]^2[u-\zeta^3(1+\zeta)^2]^2}{(u+\zeta)^2(u+\zeta^2)(u+\zeta^4)(u+1)^3}.$$
If $P=(X,Y_1)$, this gives an expression for $Y_1^\sigma=Y(\pm 2P)$.  \medskip

Since $\sigma^2 = (\zeta \rightarrow \zeta^4, u \rightarrow u)= (\zeta \rightarrow \zeta^{-1}, u \rightarrow u)$, we also have
$$Y_1^{\sigma^2}=\varepsilon^7 \frac{[u-(1+\zeta)^2][u-\zeta(1+\zeta)^2]^2[u-\zeta^2(1+\zeta)^2]^2[u-\zeta^3(1+\zeta)^2]^2}{(u+\zeta)(u+\zeta^2)(u+\zeta^3)^2(u+1)^3},$$
which coincides with $Y_2$.  We have therefore that
$$P^{\sigma^2}=(X,Y_1)^{\sigma^2}=(X,Y_2)=-P,$$
and Theorem 3.1 yields
$$P^\sigma=(X,Y_1)^\sigma = \pm 2P.$$
Since $\sigma$ is an automorphism of the extension $\mathbb{Q}(\zeta,u)/\mathbb{Q}(b)$, this shows that
$$-P=(P^\sigma)^\sigma= \pm 2P^\sigma=4P,$$
and verifies that $4P=-P$, i.e. $5P=O$.
\bigskip

\section{The Ramanujan-Rogers continued fraction.}

As in the introduction, we now set $b=r^5(\tau)$ and $\varepsilon=\frac{-1+\sqrt{5}}{2}$, where $r(\tau)$ given by (1.2) is the Rogers-Ramanujan continued fraction.  From equation (7.3) in \cite{duke} there is the identity
\begin{equation*}
r^5\left(\frac{-1}{5\tau}\right) =\frac{-r^5(\tau)+\varepsilon^5}{\varepsilon^5 r^5(\tau)+1}=\frac{-b+\varepsilon^5}{\varepsilon^5 b+1}.
\end{equation*}
Hence we have
$$r^5\left(\frac{-1}{5\tau}\right) = \frac{-b+\varepsilon^5}{\varepsilon^5 (b-\bar \varepsilon^5)}=\frac{1}{\varepsilon^5 u^5},$$
and we can take
$$u = \frac{1}{\varepsilon r\left(\frac{-1}{5\tau}\right)}.$$
On the other hand,
$$r\left(\frac{-1}{5\tau}\right)= \frac{\bar \varepsilon r(5\tau)+1}{r(5\tau)-\bar \varepsilon},$$
by (3.2) in \cite{duke}.  Hence,
$$u = \frac{r(5\tau)-\bar \varepsilon}{\varepsilon(\bar \varepsilon r(5\tau)+1)}=-\frac{r(5\tau)-\bar \varepsilon}{r(5\tau)-\varepsilon}.\eqno{(4.1)}$$
This shows that $u$ is a linear fractional expression in $r(5\tau)$, proving (1.3).  Hence the coordinates $X$ and $Y_i$ in Theorem 2.1 can be expressed as products of linear fractional expressions in $r(5\tau)$.  Since $\varepsilon$ and the coefficients of the linear fractional expressions in Theorem 2.1 lie in $\mathbb{Q}(\zeta_5)$ (see the remarks following Theorem 2.1), the same is true for $X, Y_i$ in terms of $r(5\tau)$. This also holds if $u$ is replaced by $\zeta^{-i} u$ in (4.1) while holding $\alpha=\sqrt{5}$ fixed; and letting $i$ vary yields the coordinates of $10$ of the $20$ points in $E_5(b)[5]-\langle (0,0) \rangle$.\medskip

For example, we have
$$\frac{(u-(1+\zeta)^2)}{u+1} =\frac{(1+\zeta)(1-\zeta^3)}{\sqrt{5}}(r(5\tau)-(1+\zeta^2)),$$
while
$$\frac{(u-\zeta(1+\zeta)^2)}{u+\zeta} =\zeta^2(1+\zeta) \frac{r(5\tau)-(1+\zeta)}{r(5\tau)+(1+\zeta+\zeta^2)}.$$
Further,
$$\frac{(u-\zeta^2(1+\zeta)^2)}{u+\zeta^2} = -\zeta \frac{r(5\tau)+\zeta(1+\zeta+\zeta^2)}{r(5\tau)-\zeta^2(1+\zeta^2)},$$
and
$$\frac{(u-\zeta^3(1+\zeta)^2)}{u+\zeta^3} =-\zeta(1+\zeta) \frac{r(5\tau)-(1+\zeta^3)}{r(5\tau)-\zeta(1+\zeta^2)}.$$
Using, finally, that
$$\frac{u+\zeta}{u+1}=\frac{(1-\zeta)}{\sqrt{5}}(r(5\tau)+(1+\zeta+\zeta^2)),$$
we find that the $X$-coordinate in (2.3) is given by
\begin{align*}
X&=\frac{-\varepsilon}{\sqrt{5}} \frac{[r(5\tau)-(1+\zeta)][r(5\tau)-(1+\zeta^2)][r(5\tau)-(1+\zeta^3)][r(5\tau)-(1+\zeta^4)]}{[r(5\tau)-(\zeta+\zeta^3)][r(5\tau)-(\zeta^2+\zeta^4)]}\\
&=\frac{-\varepsilon}{\sqrt{5}} \ \frac{r^4(5\tau)-3r^3(5\tau)+4r^2(5\tau)-2r(5\tau)+1}{r^2(5\tau)+r(5\tau)+\varepsilon^2},
\end{align*}
where, once again, all the ``poles'' and ``zeroes'' of this function of $r(5\tau)$ are units in $\mathbb{Q}(\zeta_5)$.  \medskip

The numerator in the last expression coincides with the numerator in Ramanujan's identity
$$
\frac{r^5(\tau)}{r(5\tau)}=\frac{r^4(5\tau)-3r^3(5\tau)+4r^2(5\tau)-2r(5\tau)+1}{r^4(5\tau)+2r^3(5\tau)+4r^2(5\tau)+3r(5\tau)+1},\eqno{(4.2)}
$$
from \cite{ber}, p. 167, while the denominator is a quadratic factor of the denominator in this identity.  (Note that $\zeta+\zeta^2, \zeta^2+\zeta^4, \zeta^3+\zeta^4, \zeta+\zeta^3$ are the conjugate roots of $x^4+2x^3+4x^2+3x+1$.) Therefore, we may also write the formula for $X$ as
$$X=\frac{-\varepsilon}{\sqrt{5}} \ \frac{r^5(\tau)}{r(5\tau)} (r^2(5\tau)+r(5\tau)+\bar \varepsilon^2).$$

The coordinates $Y_i$ in Theorem 2.1 may be computed using the above formulas, along with the formula
$$\frac{u+\zeta}{u+\zeta^4}=-\zeta \frac{r(5\tau)+1+\zeta+\zeta^2}{r(5\tau)-\zeta(1+\zeta)}.$$
We find with $\eta=(\zeta-1)/\sqrt{5}$ that
$$Y_1=\eta^3 \frac{[r(5\tau)-(1+\zeta)]^2[r(5\tau)-(1+\zeta^2)]^2[r(5\tau)-(1+\zeta^3)][r(5\tau)-(1+\zeta^4)]^2}{[r(5\tau)-(\zeta^2+\zeta^4)]^2[r(5\tau)-(\zeta+\zeta^3)][r(5\tau)-(\zeta+\zeta^2)]},$$
or
$$
Y_1=\eta^3 \frac{(r^4(5\tau)-3r^3(5\tau)+4r^2(5\tau)-2r(5\tau)+1)^2}{(r^4(5\tau)+2r^3(5\tau)+4r^2(5\tau)+3r(5\tau)+1)} \cdot Z, \eqno{(4.3)}
$$
where
$$Z=\frac{r(5\tau)-(\zeta^3+\zeta^4)}{(r(5\tau)-(\zeta^2+\zeta^4))(r(5\tau)-(1+\zeta^3))}.$$
Using (4.1) this can also be written as
$$Y_1=\eta^3 \frac{r^5(\tau)}{r(5\tau)} (r^4(5\tau)-3r^3(5\tau)+4r^2(5\tau)-2r(5\tau)+1) \cdot Z.$$
The formula for $Y_2$ can be obtained by applying the map $\sigma^2=(\zeta \rightarrow \zeta^4, u \rightarrow u)$ to $Y_1$, as in Section 3:
$$Y_2=\left(\frac{\zeta^4-1}{\sqrt{5}}\right)^3 \frac{(r^4(5\tau)-3r^3(5\tau)+4r^2(5\tau)-2r(5\tau)+1)^2}{(r^4(5\tau)+2r^3(5\tau)+4r^2(5\tau)+3r(5\tau)+1)} \cdot Z^{\sigma^2},\eqno{(4.4)}$$
with
$$Z^{\sigma^2}=\frac{r(5\tau)-(\zeta+\zeta^2)}{(r(5\tau)-(\zeta+\zeta^3))(r(5\tau)-(1+\zeta^2))}.$$
 \medskip

If we perform the same calculations by sending $\zeta$ to $\zeta^2$ in the formula (2.3), so that $\alpha=\sqrt{5}$ is replaced by $-\sqrt{5}$, $\varepsilon$ is replaced by $\bar \varepsilon$, and $u$ by $1/u$ in (4.1), then by Theorem 3.1 we find the $X$-coordinate of the double of the point $P=(X,Y_1)$:
\begin{align*}
X(2P)&=\frac{\bar \varepsilon}{\sqrt{5}} \frac{[r(5\tau)-(1+\zeta)][r(5\tau)-(1+\zeta^2)][r(5\tau)-(1+\zeta^3)][r(5\tau)-(1+\zeta^4)]}{[r(5\tau)-(\zeta^3+\zeta^4)][r(5\tau)-(\zeta+\zeta^2)]}\\
&=\frac{\bar \varepsilon}{\sqrt{5}} \ \frac{r^4(5\tau)-3r^3(5\tau)+4r^2(5\tau)-2r(5\tau)+1}{r^2(5\tau)+r(5\tau)+\bar \varepsilon^2}\\
&=\frac{\bar \varepsilon}{\sqrt{5}} \ \frac{r^5(\tau)}{r(5\tau)} (r^2(5\tau)+r(5\tau)+\varepsilon^2).
\end{align*}
The corresponding $Y$-coordinates are obtained by applying $\zeta \rightarrow \zeta^2$ to the expressions given for $Y_1, Y_2$ above. \medskip

As above, if we apply $\zeta \rightarrow \zeta^2$ to the formulas in Theorem 2.1 and set $u$ equal to the quantity
$$u=-\zeta^i \frac{r(5\tau)-\varepsilon}{r(5\tau)- \bar \varepsilon}, \ \ 0 \le i \le 4,$$
then we obtain the coordinates of the remaining $10$ points in $E_5(b)[5]-\langle (0,0) \rangle$.  This completes the proof of Theorem 1.1.

\section{Appendix.}

The roots of the equation $g(X)=0$ in Section 2 are found using Watson's method and the following quantities defined in \cite{lsw}.  First, if 
$$a_1:=\frac{1}{20}(\alpha-5)(-3-\alpha-7b+3b\alpha -2b^2)$$
is the coefficient of $X^4$ in $g(X)$, we have
$$f(x)=g(x-\frac{a_1}{5})=x^5+10Cx^3+10Dx^2+5Ex+F,$$
where:
\begin{align*}
C&=\frac{1}{4000}(-3+\alpha)(b+3+\alpha)(-4b+3+\alpha)(-2b-11+5\alpha)(2b+11+5\alpha),\\
D&=\frac{-1}{100000}(-5+2\alpha)(-8b^3-27b^2+11\alpha b^2-41b-19\alpha b+4 +4\alpha)\\
&\times (-2b-11+5\alpha)(2b+11+5\alpha)^2,\\
E&=\frac{1}{500000}(-7+3\alpha)(-6b^5-60b^4+6\alpha b^4-135b^3+47\alpha b^3+505b^2+229\alpha b^2\\
&-150b-72\alpha b+12+6\alpha)(-2b-11+5\alpha)(2b+11+5\alpha)^2,\\
F&=\frac{-1}{12500000}(-25+11\alpha)(-8b^7-133b^6+5\alpha b^6-707b^5+115\alpha b^5+3790b^4\\
&+2900\alpha b^4-15405b^3-6475\alpha b^3+5326b^2+2400 \alpha b^2-794b-360\alpha b+44+20\alpha)\\
& \times(-2b-11+5\alpha)(2b+11+5\alpha)^2.
\end{align*}
Incidentally, this shows that the polynomial $f(x)$ is irreducible over $\mathbb{Q}(\alpha,b)$, by an analogue of Eisenstein's theorem, because of the factor $-2b-11+5\alpha$ in each of the coefficients.  This yields:
\begin{align*}
K&=E+3C^2\\
&=\frac{1}{8000}(-9+4\alpha)b^2(-2b+29+13\alpha)(-2b-11+5\alpha)(2b+11+5\alpha)^2;\\
L&=-2DF+3E^2-2C^2E+8CD^2+15C^4\\
&=\frac{1}{140800000}(-35+16\alpha)b^4(-2b+1+\alpha)(22b-19+13\alpha)(-2b-11+5\alpha)^2\\
& \times (2b+11+5\alpha)^4;\\
\end{align*}
\begin{align*}
M&=CF^2-2DEF+E^3-2C^2 DF-11C^2 E^2+28CD^2E-16D^4+35C^4E\\
&-40C^3D^2-25C^6\\
&=\frac{1}{512000000}(-9+4\alpha)b^6(-2b^4-20b^3+2\alpha b^3-11b^2-11\alpha b^2-35b-13\alpha b+3+\alpha)\\
&\times (-2b-11+5\alpha)^3(2b+11+5\alpha)^5.
\end{align*}
The discriminant of $f(x)$ is
$$\delta=\frac{1}{1024000}(123-55\alpha)(-2b-11+5\alpha)^4(2b+11+5\alpha)^8 b^{14},$$
so that
$$\sqrt{\delta}=\frac{1}{1600}\alpha\left(\frac{1-\alpha}{2}\right)^5(-2b-11+5\alpha)^2(2b+11+5\alpha)^4 b^7.$$

The polynomial
$$h(x)=x^6-\frac{K}{5}x^4+\frac{L}{125}x^2-\frac{\alpha \sqrt{\delta}}{390625}x+\frac{M}{3125}$$
has the quantity
$$\theta=\frac{1}{50}b(b^2+11b-1)$$
as a root.  Treating $b$ as an indeterminate, $\theta \neq 0,\pm C$.  Hence, Theorem 1 in \cite{lsw} applies.  We take
$$T=\frac{1}{20000}(5-\alpha)b(b-2+\alpha)(-2b-11+5\alpha)(2b+11+5\alpha)^2$$
to be the solution of $p(T) = q(T) = 0$ in that theorem.  Then (2.20) in \cite{lsw} gives
\begin{align*}
R_1&=\sqrt{(D-T)^2+4(C-\theta)^2(C+\theta)}\\
&=\frac{1}{4000}(3-\alpha)b(2b-1+\alpha)(-2b-11+5\alpha)(2b+11+5\alpha)^2,
\end{align*}
and
\begin{align*}
R_2&=\frac{C(D^2-T^2)+(C^2-\theta^2)(C^2+3\theta^2-E)}{R_1\theta}\\
&=\frac{1}{2000}(-2+\alpha)b(-2b+1+\alpha)(-2b-11+5\alpha)(2b+11+5\alpha)^2.
\end{align*}
Next, we define the quantities
\begin{align*}
X'&=\frac{1}{2}(-D+T+R_1)=\frac{1}{2^65^5}(-5+2\alpha)(2b+1+\alpha)(-2b-11+5\alpha)^2(2b+11+5\alpha)^3,\\
Y&=\frac{1}{2}(-D-T+R_2)=\frac{1}{2^45^5}(-5+2\alpha)(-b+2+\alpha)^2(-2b-11+5\alpha)^2(2b+11+5\alpha)^2,\\
Z&=-C-\theta=\frac{1}{2000}(3-\alpha)(-b+2+\alpha)(2b+1+\alpha)(-2b-11+5\alpha)(2b+11+5\alpha).
\end{align*}
Then, according to Theorem 1 of \cite{lsw}, the quantity $u_1$ is the fifth root of the expression
\begin{align*}
u_1^5&=\frac{X'^2Y}{Z^2}=\frac{1}{2^{11}5^8}(-25+11\alpha)(-2b-11+5\alpha)^4(2b+11+5\alpha)^6\\
&=\frac{1}{2^{10}5^5}\left(\frac{-5+\alpha}{10}\right)^5(-2b-11+5\alpha)^5(2b+11+5\alpha)^5\\
& \times \frac{2b+11+5\alpha}{-2b-11+5\alpha}.
\end{align*}
This shows that $u_1$ is a polynomial in $b$ and $\alpha$ times $u$, where
$$u^5= \frac{2b+11+5\alpha}{-2b-11+5\alpha},$$
as in Section 2; in fact, we have
\begin{align*}
u_1 &= \frac{-5+\alpha}{200}(-2b-11+5\alpha)(2b+11+5\alpha) \times u\\
&=\frac{5-\alpha}{100}(2b^2+22b-2)u=\frac{5-\alpha}{100}A_1u.
\end{align*}
This gives the first degree term in $u$ in the expression for the root $X$ in Section 2.  Similarly, the quantities
\begin{align*}
\bar X &=\frac{1}{2}(-D+T-R_1)\\
&=\frac{1}{100000}(-5+2\alpha)(-b+2+\alpha)(-2b-1+\alpha)^2(-2b-11+5\alpha)(2b+11+5\alpha)^2,\\
\bar Y &=\frac{1}{2}(-D-T-R_2)\\
&=\frac{1}{200000}(-5+2\alpha)(-2b-1+\alpha)(2b+1+\alpha)^2(-2b-11+5\alpha)(2b+11+5\alpha)^2,\\
\bar Z &=-C+\theta=\frac{1}{4000}(3-\alpha)(-2b-1+\alpha)(-2b-11+5\alpha)(2b+11+5\alpha)^2,
\end{align*}
yield the expressions
$$u_2=\frac{\bar X}{\bar Z^2}u_1^2=\frac{5-\alpha}{100}A_2u^2, \ \ u_3=\frac{\bar X \bar Y}{Z \bar Z^3}u_1^3=\frac{5-\alpha}{100}A_3u^3,$$
$$u_4=\frac{\bar X^2 \bar Y}{Z^2 \bar Z^4}u_1^4=\frac{5-\alpha}{100}A_4u^4,$$
which are the second, third, and fourth degree terms in the expression for the root $X$, where
\begin{align*}
A_2&=(b-2-\alpha)(-2b-11+5\alpha),\\
A_3&=-\frac{1}{2}(2b+1+\alpha)(-2b-11+5\alpha),\\
A_4&=-\frac{1}{2}(-2b-1+\alpha)(-2b-11+5\alpha).
\end{align*}
Together with the fact that $\frac{(5-\alpha)}{100} A_0=-\frac{a_1}{5}$, this yields the expression for the root
$$X=u_1+u_2+u_3+u_4-\frac{a_1}{5}=\frac{(5-\alpha)}{100}(A_4u^4+A_3u^3+A_2u^2+A_1u+A_0)$$
of $g(X)=0$ in Section 2.  By replacing $u$ by $\zeta^i u$, for $0 \le i \le 4$, and solving the resulting system of linear equations for the powers of $u$, it is not hard to see that $\mathbb{Q}(\zeta,b,X)=\mathbb{Q}(\zeta,b,u)$ is the field generated over $\mathbb{Q}(\zeta,b)$ by the $X$-coordinates of the points of order $5$ on $E_5$.  This gives an alternate verification that $u^5$ is a Kummer element for the extension $\mathbb{Q}(\zeta, b, X)/\mathbb{Q}(\zeta,b)$, as in \cite{v}.

\medskip

\noindent Dept. of Mathematical Sciences, LD 270

\noindent Indiana University - Purdue University at Indianapolis (IUPUI)

\noindent Indianapolis, IN 46202

\noindent {\it e-mail: pmorton@iupui.edu}

\end{document}